\documentclass[a4paper,12pt]{article}

\usepackage[utf8]{inputenc}

\usepackage{amssymb}
\usepackage{amsmath}
\usepackage{mathbbol}
\usepackage{mathrsfs}
\usepackage{arydshln}
\usepackage{textcomp}

\usepackage[left=2cm,right=1.5cm,top=1.5cm,bottom=1.5cm,bindingoffset=0cm]{geometry}

\title{Steen-Ermakov-Pinney equation and integrable nonlinear deformation of one-dimensional Dirac equation}
\author{Y. Prykarpatskyy\\
 The Department of Applied Mathematics, Agricultural University of Krakow, Poland, \\
 Institute of Mathematics of NAS, Kyiv,	Ukraine, \\
	email: yarpry@gmail.com}

\date{}

\newtheorem{Theorem}{{\rm \bf Theorem}}[section]

\numberwithin{equation}{section}

\begin{document}

\maketitle

\begin{abstract}
The paper deals with nonlinear one-dimensional Dirac equation. We describe its invariants set 
by means of the deformed linear Dirac equation, using the fact 
that two ordinary differential equations are equivalent if their sets of invariants coincide.
\end{abstract}


\section{Introduction}

In 1874 Danish mathematician Adolph Steen wrote a paper \cite{Steen1} where
he introduced the system of two equations 
\begin{equation}
r^{\prime \prime }+qr=\frac{1}{r^{3}},\ \   \label{eqS01}
\end{equation}%
and 
\begin{equation}
\ g^{\prime \prime }+qg=0,  \label{eqS01a}
\end{equation}%
where $q=q(u)$ is a continuous function on a real interval. He discovered
that these two equations above are in some sense equivalent, that is the
general solution to the second equation (\ref{eqS01a}) gives rise to that
to the first one (\ref{eqS01}) and vise-versa. Unluckily, the paper was
published in Danish and his research was lost. Later many authors have been
rediscovering these equations and mentioning too the property above. 
In 1880 V.~Yermakov \cite{Ermakov} gave a novel derivation and
generalization of the Steen's equations. This generalization was actively
studied and developed by others researches. Later, in 1950 Edmund Pinney 
\cite{Pinney} showed that the solution of the first equation (\ref{eqS01})
is 
\begin{equation}
r(t)=\sqrt{Au^{2}+2Buv+Cv^{2}},  \label{eqS02}
\end{equation}%
where $u(t)$ and $v(t)$ are two arbitrary linearly independent solutions of
the second equation of (\ref{eqS01}), $A,B$ and $C$ are constants which
satisfy the equality 
\begin{equation}
B^{2}-AC=\frac{1}{W^{2}},  \label{eqS02a}
\end{equation}%
with $W$ being the constant Wronskian of the two independent solutions 
$u,v $ to (\ref{eqS01a}).

Nowadays the Steen's contribution is forgotten too and nobody calls systems
of equations (\ref{eqS01}) Steen's equations. The name of Ermakov is
commonly used or they also often called Ermakov-Pinney equations. Raymond
Redheffer and Irene Redheffer wrote (\cite{RRIR}) an historical survey on
the Steen's equations with English translation of the original Steen's paper.

It is worth highlighting the following quite simple fact. Let us
take the two equations 
\begin{equation}
y^{\prime \prime }=\omega (t)y  \label{eqS03}
\end{equation}%
and 
\begin{equation}
z^{\prime \prime }=\omega (t)z+\frac{k}{z^{3}},  \label{eqS04}
\end{equation}%
which are, in fact, rewritten equations from (\ref{eqS01}). These equations
can be easily transformed to the more generalized form 
$(py^{\prime})^{\prime }=\omega (t)y$ and $(pz^{\prime })^{\prime }=\omega (t)z+\frac{k}{z^{3}}$ 
for some smooth function $p(z)$ by changing the variables on which
we will not stop here.

Let us denote $\alpha :=y_{1}y_{2},$ where $y_{1}$ and $y_{2}$ are two
arbitrary solutions to the equation (\ref{eqS03}), and $\beta :=z^{2}$.
After differentiation and substituting into (\ref{eqS03}) and (\ref{eqS04})
one can obtain the following expressions: 
\begin{equation}
\alpha ^{\prime \prime \prime }-2\omega ^{\prime }\alpha -4\omega \alpha
^{\prime }=0,  \label{eqS05a}
\end{equation}%
or 
\begin{equation}
\left[ \frac{d^{3}}{dt^{3}}-\left( 2\omega \frac{d}{dt}+2\frac{d}{dt}\omega
\right) \right] \alpha =0,  \label{eqS05b}
\end{equation}%
and 
\begin{equation}
\beta ^{\prime \prime \prime }-2\omega ^{\prime }\beta -4\omega \beta
^{\prime }=0,  \label{eqS06a}
\end{equation}%
or 
\begin{equation}
\left[ \frac{d^{3}}{dt^{3}}-\left( 2\omega \frac{d}{dt}+2\frac{d}{dt}\omega
\right) \right] \beta =0.  \label{eqS06b}
\end{equation}%
It is worth mentioning that the differential expression 
$\eta =\left[ \frac{d^{3}}{dt^{3}}-\left( 2\omega \frac{d}{dt}+2\frac{d}{dt}\omega \right)\right]$ 
is the second Poisson operator \cite{Novikov,PrykMyk} for
the KdV equation $u_{t}=-u_{xxx}-6uu_{x}=K[u]$ in its Hamiltonian form 
$u_{t}=-\eta \mathrm{grad}H$, where $H$ is the corresponding \cite{Arnold}
Hamiltonian functional.

Taking into account that the equations (\ref{eqS05b}) and (\ref{eqS06b}) are
linear and the same, one can infer that their sets of the solutions coincide
to each other. It means that the a solution to the equation (\ref{eqS06b})
is the sum of the solutions to the equation (\ref{eqS05b}), whereas the
expression 
\begin{equation}
z^{2}=Ay_{1}^{2}+By_{1}y_{2}+Cy_{2}^{2}  \label{eqS07}
\end{equation}%
is the solution of the (\ref{eqS05b}), and then the function 
\begin{equation}
z=\sqrt{Ay_{1}^{2}+By_{1}y_{2}+Cy_{2}^{2}}  \label{eqS08}
\end{equation}%
is the solution of the (\ref{eqS04}). Taking into account that the Wronskian
of any two different solutions to (\ref{eqS03}) is constant, one easily
obtains the relationship (\ref{eqS02a}).

In the next section we will use the above mentioned trick in finding
solutions to the one-dimensional Dirac equations \cite{Novikov}.

\section{Dirac equation, its invariants and integrable nonlinear deformation}

Let us consider the one-dimensional Dirac equation 
\begin{equation}
\frac{df}{dx}=l(\lambda ;q)f,\ \ \ \ \ \ l(\lambda ;q):=\left( 
\begin{array}{cc}
\lambda & q_{1}(x) \\ 
q_{2}(x) & -\lambda%
\end{array}%
\right) ,  \label{eq_d01}
\end{equation}%
where $x\in \mathbb{R}$, $\lambda \in \mathbb{C}$ is a complex parameter, $%
(q_{1},q_{2})^{\intercal }\in M\simeq C^{\infty }(\mathbb{R};\mathbb{C}^{2})$
is a functional vector of potentials, and $f\in (L_{\infty }(\mathbb{R};%
\mathbb{C}^{2})$ is found solution to (\ref{eq_d01}). The whole solution set
to the equation (\ref{eq_d01}) is completely described \cite{CoLe} by means
of the corresponding fundamental solution 
$F:=\left( 
\begin{array}{cc}
f_{11} & f_{12} \\ 
f_{21} & f_{22}%
\end{array}%
\right) \in C^{\infty }(\mathbb{R}\times \mathbb{R},\mathrm{End}\mathbb{C}%
^{2})$, 
satisfying the matrix equation 
\begin{equation}
\frac{d}{dx}F(x,x_{0})=l(\lambda ;q)F(x,x_{0}),\ \ \ \ \ F(x,x_{0})\big|%
_{x=x_{0}}=\mathbf{1}  \label{eq_d02}
\end{equation}%
at any point $x_{0}\in \mathbb{R}.$ Then an arbitrary solution $f\in
C^{\infty }(\mathbb{R};\mathbb{C}^{2})$ to (\ref{eq_d02}), evidently, can be
represented as 
\begin{equation}
f(x)=F(x,x_{0})f(x_{0}),  \label{eq_d03}
\end{equation}%
where $f(x_{0})\in \mathbb{C}^{2}$ is a suitable Cauchy data vector.

It is well known from the general theory of ordinary differential equations 
\cite{CoLe} that the solution set to the equation (\ref{eq_d01}) can be
equivalently described by means of its complete set of invariants. Moreover,
the two ordinary differential equations are then considered to be equivalent
if their sets of invariants coincide. From this point of view one can
consider a suitable deformed Dirac equation 
\begin{equation}
\frac{d\tilde{f}}{dx}=l(\lambda ;q)\tilde{f}+\delta \tilde{f},
\label{eq_d04}
\end{equation}%
where a vector $\delta \tilde{f}\in C^{\infty }(\mathbb{R}\mathbb{C}^{2})$
can depend on $\tilde{f}\in C^{\infty }(\mathbb{R};\mathbb{C}^{2})$ and the
vector $q\in M$ coincides with that chosen in \ (\ref{eq_d01}) .

Now a problem consists in determining the vector $\delta \tilde{f}\in
C^{\infty }(\mathbb{R};\mathbb{C}^{2})$ in such a form which guarantees that
the set of invariants of (\ref{eq_d04}) will contain or coincide with that
of the equation (\ref{eq_d01}).

To approach a solution to this problem we need first to describe the
invariants set to the equation (\ref{eq_d01}). To do that we assume for
simplicity that the functional vector $q\in M$ is $2\pi $-periodic in $x\in 
\mathbb{R}$: $q(x+2\pi )=q(x)$ for all $x\in \mathbb{R}$. Then one can
define \cite{Titch} the monodromy matrix $S(x)\in \mathrm{End}\mathbb{C}^{2}$
as 
\begin{equation}
S(x):=F(x+2\pi ,x),  \label{eq_d05}
\end{equation}%
satisfying the well known Novikov commutator equation \cite{Nov} 
\begin{equation}
\frac{dS(x)}{dx}=[l(\lambda ;q),S(x)].  \label{eq_d06}
\end{equation}%
As a consequence from (\ref{eq_d06}) one easily obtains that all functions 
\begin{equation}
\gamma _{j}:=\mathrm{tr}S^{j}(x),\ \ \ \ \frac{d\gamma _{j}}{dx}=0,
\label{eq_d07}
\end{equation}%
where $j\in \mathbb{Z}_{+}$, are invariants for (\ref{eq_d01}) and form 
\cite{Reiman} its complete set. In particular, we can determine only two
dependent invariants for (\ref{eq_d01}) 
\begin{equation}
\gamma _{1}=\mathrm{tr}S(x),\ \ \ \ \gamma _{2}=\det S(x)=1,  \label{eq_d08}
\end{equation}%
and will now try to search for such a deformation vector $\delta \tilde{f}%
\in C^{\infty }(\mathbb{R};\mathbb{C}^{2})$ which will give rise to the
invariants set of (\ref{eq_d04}) coinciding with (\ref{eq_d08}) of (\ref%
{eq_d01}). For this problem to be solved effectively, we need to find the
determining equations for invariants (\ref{eq_d08}) as for the independent
gradient $\mathrm{grad}\gamma _{1}\in T^{\ast }(M)$, 
\begin{equation}
\left( 
\begin{array}{cc}
q_{2}D_{x}^{-1}q_{2} & \frac{d}{dx}+\lambda -q_{2}D_{x}^{-1}q_{1} \\ 
\frac{d}{dx}-\lambda -q_{1}D_{x}^{-1}q_{2} & q_{1}D_{x}^{-1}q_{1}%
\end{array}%
\right) \mathrm{grad}\gamma _{1}=0,  \label{eq_d09}
\end{equation}%
where, by definition, we put $D_{x}^{-1}(\cdot ):=\frac{1}{2}\left(
\int\limits_{x_{0}}^{x}(\cdot )dy-\int\limits_{x}^{x_{0}+2\pi }(\cdot
)dy\right) ,x\in \mathbb{R}.$

It is now worth observing that the monodromy matrix $S(x)\in \mathrm{End}%
\mathbb{C}^{2}$ for any $x\in \mathbb{R}$ allows the following matrix
representation: 
\begin{equation}
S(x)=F(x,x_{0})C(x_{0})F^{-1}(x,x_{0})  \label{eq_d10}
\end{equation}%
for some specially chosen matrix $C(x_{0})\in \mathrm{End}\mathbb{C}^{2}$, $%
x_{0}\in \mathbb{R}.$ As a simple consequence of the representation (\ref%
{eq_d10}) one easily obtains that 
\begin{equation}
\mathrm{grad}\gamma _{1}=\left( 
\begin{array}{l}
c_{11}f_{21}f_{22}-c_{12}f_{21}^{2}+c_{21}f_{22}^{2}-c_{22}f_{22}f_{21} \\ 
c_{12}f_{11}^{2}-c_{11}f_{11}f_{12}-c_{21}f_{12}^{2}+c_{22}f_{12}f_{11}%
\end{array}%
\right) ,  \label{eq_d11}
\end{equation}%
where we used the matrix expression $C(x_{0})=\left( 
\begin{array}{cc}
c_{11} & c_{12} \\ 
c_{21} & c_{22}%
\end{array}%
\right) \in \mathrm{End}\mathbb{C}^{2}$. Taking into account the
arbitrariness of the matrix $C(x_{0})\in \mathrm{End}\mathbb{C}^{2}$
entering (\ref{eq_d10}), we can easily obtain from (\ref{eq_d09}) and (\ref%
{eq_d11}) that the following equation 
\begin{equation}
\left( 
\begin{array}{cc}
q_{2}D_{x}^{-1}q_{2} & \frac{d}{dx}+\lambda -q_{2}D_{x}^{-1}q_{1} \\ 
\frac{d}{dx}-\lambda -q_{1}D_{x}^{-1}q_{2} & q_{1}D_{x}^{-1}q_{1}%
\end{array}%
\right) \left( 
\begin{array}{c}
-f_{21}^{2} \\ 
f_{11}^{2}%
\end{array}%
\right) =0  \label{eq_d12}
\end{equation}%
holds for all $x\in \mathbb{R}$ and $\lambda \in \mathbb{C}$.

Now it is easy to infer that if the vector 
$\tilde{f}:=(\tilde{f}_{1},\tilde{f}_{2})^{\tau }\in ^{\infty }(\mathbb{C};\mathbb{C}^{2})$ 
of the equation (\ref{eq_d04}) satisfies the same equation as (\ref{eq_d12}) 
\begin{equation}
\left( 
\begin{array}{cc}
q_{2}D_{x}^{-1}q_{2} & \frac{d}{dx}+\lambda -q_{2}D_{x}^{-1}q_{1} \\ 
\frac{d}{dx}-\lambda -q_{1}D_{x}^{-1}q_{2} & q_{1}D_{x}^{-1}q_{1}%
\end{array}%
\right) \left( 
\begin{array}{c}
-\tilde{f}_{2}^{2} \\ 
\tilde{f}_{1}^{2}%
\end{array}%
\right) =0,  \label{eq_d13}
\end{equation}%
then the corresponding set of invariants of the deformed equation 
(\ref{eq_d04}) will possess that of invariants for (\ref{eq_d01}). In particular,
from (\ref{eq_d13}) it follows that a partial solution to the deformed Dirac
equation (\ref{eq_d04}) can be represented as 
\begin{equation}
\tilde{f}=\left( 
\begin{array}{l}
(c_{12}f_{11}^{2}-c_{21}f_{11}f_{12}-c_{21}f_{12}^{2}+c_{22}f_{12}f_{11})^{1/2}
\\ 
(c_{12}f_{21}^{2}-c_{11}f_{21}f_{22}+c_{22}f_{22}f_{21}-c_{21}f_{22}^{2})^{1/2}%
\end{array}%
\right) ,  \label{eq_d14}
\end{equation}%
depending only on the fundamental matrix $F(x,x_{0})\in \mathrm{End}\mathbb{C%
}^{2}$ of the equation (\ref{eq_d01}) and equivalently, on its set of
invariants. It is clear that the deformed Dirac equation (\ref{eq_d04}) can
generate a new solutions to it, yet the problem of describing this set of
invariants is much more complicated and will not be herewith discussed.

Let us proceed now to describing the deformed Dirac equation (\ref{eq_d04}),
taking into account that the vector $(-\tilde{f}_{2}^{2},\tilde{f}%
_{1}^{2})\in C^{\infty }(\mathbb{R};\mathbb{C}^{2})$ satisfies in general
the following equation: 
\begin{equation}
\left( 
\begin{array}{cc}
q_{2}D_{x}^{-1}q_{2} & \frac{d}{dx}+\lambda -q_{2}D_{x}^{-1}q_{1} \\ 
\frac{d}{dx}-\lambda -q_{1}D_{x}^{-1}q_{2} & q_{1}D_{x}^{-1}q_{1}%
\end{array}%
\right) \left( 
\begin{array}{c}
-\tilde{f}_{2}^{2} \\ 
\tilde{f}_{1}^{2}%
\end{array}%
\right) =\left( 
\begin{array}{c}
-\delta \tilde{f}_{2}\ \tilde{f}_{2}+q_{1}\frac{d^{-1}}{dx}(\delta \tilde{f}%
_{1}\ \tilde{f}_{2}+\delta \tilde{f}_{2}\ \tilde{f}_{1}) \\ 
\delta \tilde{f}_{1}\ \tilde{f}_{1}-q_{2}\frac{d^{-1}}{dx}(\delta \tilde{f}%
_{1}\ \tilde{f}_{2}+\delta \tilde{f}_{2}\ \tilde{f}_{1})%
\end{array}%
\right) ,  \label{eq_d15}
\end{equation}%
which reduces to (\ref{eq_d13}) if identically one has 
\begin{eqnarray}
\delta \tilde{f}_{2}\ \tilde{f}_{2} &=&q_{1}D_{x}^{-1}(\delta \tilde{f}_{1}\ 
\tilde{f}_{2}+\delta \tilde{f}_{2}\ \tilde{f}_{1}),  \label{eq_d16} \\
\delta \tilde{f}_{1}\ \tilde{f}_{1} &=&q_{2}D_{x}^{-1}(\delta \tilde{f}_{1}\ 
\tilde{f}_{2}+\delta \tilde{f}_{2}\ \tilde{f}_{1}).  \notag
\end{eqnarray}%
As a simple consequence of (\ref{eq_d16}) one obtains that there exists some
function $\alpha \in C^{\infty }(\mathbb{R};\mathbb{C})$ for which 
\begin{equation}
\delta \tilde{f}_{2}=\frac{\alpha }{\tilde{f}_{2}}q_{1},\ \ \ \ \delta 
\tilde{f}_{1}=\frac{\alpha }{\tilde{f}_{1}}q_{2}.  \label{eq_d17}
\end{equation}%
Having substituted (\ref{eq_d17}) into (\ref{eq_d16}) one obtains that 
\begin{equation}
\frac{d\alpha }{dx}=\alpha \left( q_{1}\frac{\tilde{f}_{1}}{\tilde{f}_{2}}%
+q_{2}\frac{\tilde{f}_{2}}{\tilde{f}_{1}}\right) ,  \label{eq_d18}
\end{equation}%
or upon integrating (\ref{eq_d18}), one ensues 
\begin{equation}
\alpha =\bar{\alpha}\exp [D_{x}^{-1}(q_{1}\frac{\tilde{f}_{1}}{\tilde{f}_{2}}%
+q_{2}\frac{\tilde{f}_{2}}{\tilde{f}_{1}})],  \label{eq_d19}
\end{equation}%
where $\bar{\alpha}\in \mathbb{C}$ is an arbitrary constant. Having now
summarized the results obtained above one can formulate the following
theorem.

\begin{Theorem}\label{th01} 
	Consider two Dirac type equations: the first one (\ref{eq_d01})
	linear and the second one (\ref{eq_d04}) nonlinear, where 
	\begin{eqnarray}
	\delta \tilde{f}_{1} &=&\frac{\bar{\alpha}}{\tilde{f}_{1}}q_{2}\exp \left[
	D_{x}^{-1}(q_{1}\frac{\tilde{f}_{1}}{\tilde{f}_{2}}+q_{2}\frac{\tilde{f}_{2}%
	}{\tilde{f}_{1}})\right] ,  \label{eq_d20} \\
	&&  \notag \\
	\delta \tilde{f}_{2} &=&\frac{\bar{\alpha}}{\tilde{f}_{2}}q_{1}\exp \left[
	D_{x}^{-1}(q_{1}\frac{\tilde{f}_{1}}{\tilde{f}_{2}}+q_{2}\frac{\tilde{f}_{2}%
	}{\tilde{f}_{1}})\right]  \notag
	\end{eqnarray}%
	and $\bar{\alpha}\in \mathbb{C}$ is an arbitrary constant. Then a partial
	solution to the nonlinear Dirac type equation (\ref{eq_d04}) is given by the
	explicit expression (\ref{eq_d14}), represented by means of the fundamental
	solution $F(x,x_{0})\in \mathrm{End}\mathbb{C}^{2}$ to the Dirac equation (%
	\ref{eq_d01}) and the arbitrary constant matrix $C\in \mathrm{End}\mathbb{C}%
	^{2}$.
\end{Theorem}

Thus, the Dirac equation (\ref{eq_d04}), deformed by means of the vector
components (\ref{eq_d20}), is a nonlinear integro-differential equation
depending on the functional element $q\in M$, whose $2\pi $-periodicity
assumed before is not essential, as the main inferences, which are presented
above, were based strictly on local reasonings.

\section{Conclusion}

Based on the analogy with the oscillator type equations (\ref{eq_d01}) and (\ref{eq_d02}) we
succeeded in deriving a more general Steen type statement about the
relationship between the solution sets to the linear Dirac equation 
(\ref{eq_d01}) and its nonlinear deformation (\ref{eq_d04}), specified by the
expressions (\ref{eq_d20}).


\begin{thebibliography}{99}
	
\bibitem{Steen1} A. Steen, Om Formen for Integralet af den lineaere Differentialligning af an den Orden. Overs. over d. K. Danske Vidensk. Selsk. Forh. 1874, pp. 1-12

\bibitem{Ermakov} V. Ermakov, Second-order differential equations. Conditions of complete integrability,
Universitetskie Izvestiya, Kiev, 1880, N 9, 1-25 (translated by A.O. Harin)

\bibitem{Pinney} E. Pinney, The nonlinear differential equation $y''(x)+p(x)y+ cy^{-3}= 0$, Proc. Amer. Math. Soc., 1950, V.1, 681.

\bibitem{RRIR} R. Redheffer, I. Redheffer, Steen's 1847 paper: historical survey and translation. Aequationes Math. 2001, 61, 2001, pp. 131-150

\bibitem{Novikov} S.P. Novikov (Ed.), Theory of Solitons, Moscow, Nauka, 1980


\bibitem{PrykMyk} A. Prykarpatsky, I. Mykytyuk, Algebraic Integrability of Nonlinear Dynamical Systems on Manifolds: Classical and Quantum Aspects, Kluwer Academic Publishers, the Netherlands, 1998. 




\bibitem{Arnold} V.I. Arnold, Mathematical methods of classical mechanics, Springer, New York, 1978. 
	
\bibitem{CoLe}  E.A. Coddington,N. Levinson, Theory of Ordinary Differential
Equations, McGraw-Hill, New York, 1955.

\bibitem{Titch} E.C. Titchmarsh, Eigenfunction  expansions  associated  with  second-order  differential
equations, Part One. Second edition. Oxford University Press, Oxford, 1962.

\bibitem{Nov} S. Novikov, S.V. Manakov, L.P. Pitaevskii, V.E. Zakharov, Theory of solitons. The inverse scattering methods.
Monographs in Contemporary Mathematics. Springer, 1984



\bibitem{Reiman} A.G. Reiman, M.A. Semenov-Tyan-Shanskii, Current algebras and nonlinear partial differential equations, Dokl. Akad. Nauk SSSR,251, No. 6, 1980, pp. 1310-1313 


\end{thebibliography}
\end{document}